\RequirePackage{fix-cm}
\documentclass[opre,sglanonrev]{informs4}
\RequirePackage{newtxtext}
\RequirePackage{newtxmath}
\usepackage{tabularx}
\OneAndAHalfSpacedXI

\usepackage[ruled]{algorithm2e}
\usepackage{tikz}

\usepackage{natbib}
 \bibpunct[, ]{(}{)}{,}{a}{}{,}%
 \def\bibfont{\small}%

\def\OPT{{\mathrm{OPT}}}
\newcommand{\cN}{\mathbf{N}}
\newcommand{\cC}{\mathcal{C}}
\newcommand{\vecc}{\mathbf{c}}
\newcommand{\vecu}{\mathbf{u}}
\newcommand{\vecv}{\mathbf{v}}
\newcommand{\vecw}{\mathbf{w}}
\newcommand{\bfmu}{\boldsymbol{\mu}}
\newcommand{\E}{\mathbb{E}}

\EquationsNumberedThrough

\TheoremsNumberedThrough
\MANUSCRIPTNO{}
\def\theARTICLETOPLEFT{}\def\theARTICLETOPRIGHT{}

\begin{document}

\RUNAUTHOR{Sun and Udwani}

\RUNTITLE{Dynamic Assortment Optimization Under Arbitrary Demand Distributions}

\TITLE{Dynamic Assortment Optimization Under Arbitrary Demand Distributions: The Uniform Price Setting}

\ARTICLEAUTHORS{%
\AUTHOR{Shuo Sun}
\AFF{Department of Industrial and Systems Engineering,
University of Wisconsin--Madison, Madison, WI,
\EMAIL{ssun279@wisc.edu}}

\AUTHOR{Rajan Udwani}
\AFF{Department of Industrial Engineering and Operations Research,
University of California, Berkeley, Berkeley, CA 94720,
\EMAIL{udwani@berkeley.edu}}
}
\ABSTRACT{%
We study the joint optimization of assortment and inventory when products are identically priced and the total demand (number of customers) is uncertain.
A seller chooses initial inventories before customers arrive and every customer buys her preferred product among those still in stock according to a choice model known to the seller. The objective is to maximize the (expected) total sales, subject to a capacity constraint. We give a $(1-1/e)^2$-approximation under any ranking-based choice model and a $(1-1/e-\varepsilon)$-approximation under the Markov chain choice model, where $\varepsilon$ is a tunable parameter that trades off runtime with approximation guarantee. Our guarantees hold for arbitrary demand distributions and do not require the ubiquitous Increasing Failure Rate (IFR) assumption. In fact, our results also improve on the state-of-the-art approximations for IFR demand distributions. Key to our results is the discovery of hidden submodularity in these problems and the development of efficient procedures to evaluate (or approximate) the objective values. Specifically, we uncover a novel connection to online bipartite matching (for ranking-based models) and use stochastic coupling arguments (for the Markov chain choice model) to reduce the problems to instances of constrained submodular maximization. This opens the door to applying the vast algorithmic machinery for submodular maximization to our problem setting, leading, for example, to natural generalizations beyond cardinality constraint where no approximation results were previously known.
}%

\maketitle

\section{Introduction}
Assortment optimization is about finding the set of products a seller should offer to maximize revenue when customers substitute among them according to a discrete choice model.
When products have different prices, theoretically, the problem is very challenging to even approximate: under ranking-based choice models (defined below), it is NP-hard to approximate within a factor of $n^{1-\varepsilon}$ for $n$ products and any fixed $\varepsilon>0$~\citep{aouad2018approximability}. A large body of work is dedicated to the development of provably good efficient (polynomial time) algorithms for various (more structured) families of choice models~\citep{talluri2004revenue,rusmevichientong2010dynamic,davis2014assortment,blanchet2016markov,feldman2017revenue,desir2020constrained,berbeglia2020assortment,desir2014near}.

When product prices are uniform (identical),\footnote{Prices are often identical when products are horizontally differentiated and differ mainly in attributes like styles and sizes of apparel, colors or configurations of consumer electronics, and cosmetics shades.} or alternatively, when the objective is to maximize sales, the problem reduces to an instance of submodular maximization for a very general family of ranking-based choice models~\citep{berbeglia2020assortment}. In these models, each customer draws a random preference ranking and buys her most preferred offered product, and the purchase probability is a monotone submodular function of the assortment. This unlocks the vast number of algorithms available for submodular optimization and has enabled the development of efficient algorithms with constant-factor approximation guarantees for a variety of practical generalizations, such as ranking products under group fairness constraints~\citep{asadpour2023sequential} and product placement~\citep{elhousni2025placement}.

While the classic \emph{static} formulation
assumes that inventory is easily replenished, or alternatively that a single customer chooses from the offered assortment, we are interested in the \emph{dynamic} setting, where inventory is limited and is not replenished during the selling season. Before the season, the seller chooses an initial inventory vector $\vecc=(c_1,\dots,c_n)$ for $n$ products, subject to a capacity constraint $\sum_{i:c_i>0}c_i\le K$. Then a random number $T$ of customers arrive one at a time and the distribution of $T$ is known to the seller. Each customer independently draws her preferences from a known choice model and buys her most preferred product among those still in stock, or leaves without buying. The goal is to maximize the expected number of sales $F(\vecc)$ and we let $\OPT$ denote the optimal value. Each purchase depletes the stock, and a customer whose preferred product has sold out may substitute to another product or leave. So the assortment faced by each customer is determined by the purchases of earlier customers, and the value of a unit depends on how it changes the substitution of later customers. This \emph{stockout-based substitution}, studied since \citet{mahajan2001stocking}, makes the problem much harder than its static counterpart~\citep{goyal2016near,aouad2018greedy,aouad2019approximation,GoldsteinFeder2023,sun2024unified,fu2026joint,akchen2026inventory,zuo2026oneshot}.

 However, unlike the static case, no similar solution paradigm that works for a general choice model, for arbitrary demand distributions, and for many constraint families is known in the uniform price dynamic setting, and many fundamental questions remain open.%
 \footnote{When $T$ is sampled from a distribution with the Increasing Failure Rate (IFR) property, \citet{fu2026joint} show that the problem can be ``reduced'' (with some loss) to capacity-constrained static assortment optimization under general substitutable choice models (see Section~\ref{sec:related}). However, IFR excludes many high-variance distributions that may arise in practice~\citep{bai2023fluid,aouad2022nonparametric}.}  We believe this is due to two main obstacles:
\begin{enumerate}
\item \emph{Absence of submodular-like structure.} In the static setting with uniform prices, submodularity is the structure that enables various generalizations and extensions. Under stockout-based substitution, expected sales need not be submodular in the inventory vector under a ranking-based model~\citep[Lemma~A.4]{aouad2019approximation}, even with just two rankings (we revisit this in Section~\ref{sec:rcm}).
 In fact, even the deterministic fluid version of stockout-based substitution, in which customers with each ranking arrive as a continuous flow and buy their most preferred product in stock, is not submodular. Extensions of submodularity, such as submodular order~\citep{udwani2025submodular}, do not apply either.
 As a result, prior work has developed results tailored to specific choice models and constraint families, under additional restrictions on the distribution of $T$.
\item \emph{Evaluation.} Even when the objective is submodular, achieving the optimal guarantee $1-1/e$ of submodular maximization~\citep{nemhauser1978} requires estimates of the objective with sufficiently small error~\citep{streeter2008online,horel2016maximization}. Under the multinomial logit (MNL) model, expected sales are submodular for every distribution of $T$~\citep[Section~3.2.4]{aouad2018greedy}. However, the best known polynomial-time guarantee under MNL choice with random $T$ is strictly lower than $1-1/e$ and based on more sophisticated algorithms, because of the lack of an efficient way to evaluate the expected sales. Efficient evaluation has been noted as an open problem even for simple models such as interval preference lists~\citep{aouad2019approximation}. In fact, some results are conditioned on access to an oracle for the objective~\citep{goyal2016near,aouad2018greedy}.

Recently, \citet{zhou2026optimal} showed that, under MNL, stockout-based substitution is equivalent to a stochastic process driven by independent Poisson clocks, which turns expected sales into a one-dimensional integral that, we believe, could be evaluated to any desired accuracy. However, this relies heavily on the independence of irrelevant alternatives (IIA) property of MNL and does not extend to other choice models.
\end{enumerate}

\subsection{Our Results}\label{sec:results}
We overcome both obstacles using remarkably simple new ideas. For ranking-based models, we make a new connection with the seemingly unrelated online bipartite matching problem. Using this connection we propose a simple LP-based surrogate objective that
is monotone and submodular in the inventory vector, and
approximates the true objective within a factor of $1-1/e$ for every number of customers. Maximizing the surrogate with the standard greedy algorithm for monotone submodular functions~\citep{nemhauser1978} gives the following result.

\begin{theorem}\label{thm:rcm-main}
Under any ranking-based choice model and for every distribution of $T$, there is a deterministic algorithm that returns an inventory vector $\hat\vecc$ with $\sum_i\hat c_i\le K$ and $F(\hat\vecc)\ge(1-1/e)^2\,\OPT$ in time polynomial in $n$, $K$ and the size of the support of $T$.
\end{theorem}
We note that the ranking-based choice model may be given as an explicit list of polynomially many rankings, or through an oracle for the no-purchase probabilities $\phi(0,S)$. In the latter case, our algorithm makes only a polynomial number of calls to the oracle.

The factor $(1-1/e)^2$ combines the factor $1-1/e$ of greedy for maximizing a monotone submodular function with the factor $1-1/e$ lost by approximating the objective with a monotone submodular surrogate. Previously, for ranking-based models with uniform prices, \citet{aouad2019approximation} gave a $0.193$-approximation when $T$ is IFR, which \citet{GoldsteinFeder2023} improved to $0.226$. Very recently (concurrent with our work), \citet{akchen2026inventory} obtain a $0.4328(1-1/e)\approx0.274$-approximation when $T$ is IFR. Theorem~\ref{thm:rcm-main} improves on these guarantees, since $(1-1/e)^2\approx0.399$, and it gives the first constant-factor guarantee for ranking-based models when the distribution of $T$ is arbitrary. On the other hand, no polynomial-time algorithm can improve on $1-1/e$ unless $\mathrm P=\mathrm{NP}$, even when $T\equiv1$~\citep{goyal2016near}, because the problem then contains maximum coverage~\citep{feige1998threshold}.

Given that expected sales are submodular under MNL but not under general ranking-based models, it is natural to ask which more general models retain submodularity of the original objective. By generalizing the stochastic coupling approach for MNL~\citep{aouad2018greedy}, we show that for every distribution of $T$, the sales objective is monotone submodular for the Markov chain model~\citep{blanchet2016markov}. However, recall that submodularity by itself does not give an improved guarantee, because algorithms for submodular maximization need sufficiently accurate estimates of the objective. We also develop a general estimation procedure and this leads to the following improved result for the Markov chain model, and hence for MNL.

\begin{theorem}\label{thm:greedy}
Under any Markov chain choice model, for every distribution of $T$ and every $\varepsilon,\delta\in(0,1)$, there is a randomized algorithm that returns an inventory vector $\hat\vecc$ with $\sum_i\hat c_i\le K$ and, with probability at least $1-\delta$, $F(\hat\vecc)\ge(1-1/e-\varepsilon)\OPT$. Its running time is polynomial in $n$, $K$, $1/\varepsilon$, $\log(1/\delta)$, and the size of the support of $T$.
\end{theorem}
Given submodularity, Theorem~\ref{thm:greedy} follows from the guarantee of the greedy algorithm when objective values are only estimated~\citep[Theorem~3]{streeter2008online}, applied with our estimates of expected sales. The parameter $\varepsilon$ is the error in estimating the objective. A smaller $\varepsilon$ requires more samples in the Monte Carlo estimation and hence more running time.

Since $1-1/e$ is the best known guarantee for the Markov chain model even in the static (single customer) case~\citep{desir2020constrained}, our result indicates that the dynamic problem with an arbitrary distribution of $T$ \emph{may not be any harder than the static problem}. Previously, the best known guarantee under the Markov chain model was the same as for ranking-based models, which is $\approx0.274$~\citep{akchen2026inventory} when $T$ is IFR.

\begin{remark}[Running Time Dependence on the Support of $T$]\label{rem:support}
	Given the focused nature of this paper, we let the running times in Theorems~\ref{thm:rcm-main} and~\ref{thm:greedy} depend on the distribution of $T$ through the size of its support. For ranking-based models, the runtime dependence on $T$ can (essentially) be removed without changing the guarantee (see Remark~\ref{rem:parametric}). For the Markov chain model, by losing an $\varepsilon$ in the guarantee, we can make the runtime polynomial in $\varepsilon^{-1}\log t_{\max}$, where $t_{\max}$ is the largest value in the support of $T$ (see Remark~\ref{rem:horizon}). For parametric distributions with infinite support, such as Poisson, no truncation is needed for ranking-based models. For the Markov chain model, one can again lose an $\varepsilon$ in the guarantee to truncate the support at a sufficiently large $t_{\max}$.
\end{remark}

\begin{remark}[Running Time Dependence on $K$]\label{rem:runtime}
The running time in Theorem~\ref{thm:rcm-main} can be made polynomial in $\log K$ instead of $K$ by losing an $\varepsilon$ in the guarantee. Specifically, we can replace greedy with the algorithm of \citet[Theorem~1]{soma2018maximizing} for maximizing a monotone DR-submodular function. This gives a $((1-1/e)^2-\varepsilon)$-approximation in time polynomial in $n$, $\log K$, $1/\varepsilon$ and the size of the support of $T$. However, this does not directly apply to Theorem~\ref{thm:greedy}, where expected sales are estimated by sampling. Both standard greedy and the algorithm of \citet{soma2018maximizing} compare the marginal gains of single units, so an error of $\eta$ in each estimate can cost up to $K\eta$ in the guarantee. Therefore, the estimates need to be accurate up to about $\varepsilon\OPT/K$, which requires a number of samples that is polynomial in $K$.
\end{remark}

\begin{remark}[Budget and Matroid Constraints]\label{rem:constraints}
Our approach also applies to a budget constraint $\sum_iw_ic_i\le W$ and to matroid constraints on units of inventory (that arise in product placement, for example,~\citep{elhousni2025placement}). If we represent each unit as a separate element, our objectives are monotone submodular set functions. Therefore, we can apply partial enumeration for a budget constraint~\citep{sviridenko2004} and continuous greedy with pipage rounding for a matroid constraint~\citep{calinescu2011}. Both algorithms tolerate small errors in value queries~\citep{golovin2014online,nie2023framework}. This gives the guarantees $(1-1/e)^2-\varepsilon$ and $1-1/e-\varepsilon$, with high probability for the randomized algorithms.
To illustrate the main ideas, we focus on the capacity constraint in the rest of the paper.
\end{remark}

\subsection{Related Work}\label{sec:related}
Given the focused nature of this paper, we discuss only the most closely related work and refer to the surveys of \citet{kok2015assortment} and \citet{strauss2018review} and to the papers cited below for broader discussions. We first discuss results for the dynamic problem with a capacity constraint, grouped by the assumption on the distribution of $T$, and then results for other constraints.

\paragraph{Arbitrary Distributions of $T$.}
Obtaining guarantees for general distributions of $T$ was posed as an open question by \citet{aouad2018greedy}. In applications such as retail, the number of customers may have high variance~\citep{bai2023fluid,aouad2022nonparametric} and need not be IFR. Without distributional assumptions, \citet{aouad2019approximation} give a $(1-1/e)$-approximation for nested preference lists and \citet{segev2019assortment} gives a quasi-polynomial-time approximation scheme. Under MNL with general prices, \citet{aouad2022stability} give an approximation scheme whose running time is polynomial in all parameters except for an exponential dependence on the logarithm of the ratio between the largest and smallest MNL preference weights. When the seller can also personalize the assortment offered to each customer, \citet{sun2024unified} give a $(0.25-\varepsilon)$-approximation under a mixture of polynomially many MNL models for arbitrary distributions of $T$, relative to an upper bound on the optimal personalized policy. With uniform prices, their personalization policy offers every available product to every customer, so this guarantee applies to our problem. Therefore, prior to our work, constant-factor guarantees for arbitrary distributions of $T$ were known only for nested preference lists and mixtures of polynomially many MNL models.

\paragraph{IFR Distributions of $T$.}
Most approximation guarantees for the dynamic problem assume that $T$ has an IFR distribution, a class that includes deterministic and Poisson $T$. For nested preference lists, \citet{goyal2016near} give a polynomial-time approximation scheme. For ranking-based models, \citet{aouad2019approximation} give the $0.193$-approximation for uniform prices discussed above, which \citet{GoldsteinFeder2023} improve to $0.226$, and an $O(\log(p_n/p_1))$-approximation for general prices $p_1\le\dots\le p_n$. Under MNL with general prices, \citet{aouad2018greedy} give a $(0.139-\varepsilon)$-approximation assuming an oracle that evaluates expected revenue, and an implementable $(0.122-\varepsilon)$-approximation without it. \citet{sun2024unified} improve the guarantee to $0.286-\varepsilon$. Subsequently, \citet{zuo2026oneshot} improve it to $1/2-\varepsilon$.

Beyond MNL, \citet{fu2026joint} reduce the problem with general prices to cardinality-constrained static assortment optimization under general substitutable choice models. An $\alpha$-approximation for the static problem gives about $0.25\alpha$. Recently, \citet{feldman2026inventory} improve the factor to $0.4328\alpha$, which gives $0.4328(1-1/e)\approx0.274$ for ranking-based models with uniform prices, where greedy gives $\alpha=1-1/e$.

With uniform prices, Theorems~\ref{thm:rcm-main} and~\ref{thm:greedy} improve on all of these IFR guarantees, except the approximation scheme for nested preference lists.

\paragraph{Deterministic $T$ and Other Special Cases.}
When $T$ is deterministic, \citet{sun2024unified} obtain $0.474-\varepsilon$ under MNL with general prices, which \citet{zuo2026joint} improve to $1/2-\varepsilon$ and, subsequently, \citet{zuo2026oneshot} to $1-1/e-\varepsilon$.\footnote{With uniform prices, this guarantee essentially follows from \citet{sun2024unified}. In this case, the fluid depletion of any inventory vector attains the value of their fluid relaxation. Their comparison of stochastic and fluid sales then shows that expected sales are at least $1-1/e$ times this value, and their algorithm approximately maximizes this value over inventory vectors.} When $T$ is deterministic or Poisson, the reduction of \citet{fu2026joint} gives about $0.272\alpha$, which \citet{akchen2026inventory} improve to $0.486\alpha$, about $0.307$ for ranking-based models with uniform prices. Under MNL with identical prices and preference weights, a balanced inventory vector is optimal when $T$ is deterministic~\citep{zhang2024leveraging,zhou2026optimal}.

\paragraph{Beyond a Capacity Constraint.}
Under MNL, \citet{sun2024unified} obtain the same guarantees under a budget constraint, resolving an open problem of \citet{aouad2018greedy}. For allocating inventory across a network of stores with product-supply and store-capacity constraints, under MNL with deterministic $T$, \citet{zuo2026joint} give a $((1-1/e)/2-\varepsilon)$-approximation when each store commits upfront to the products it displays, which \citet{zuo2026oneshot} improve to $(1-1/e)^2-\varepsilon$.

\paragraph{Outline.}
In Section~\ref{sec:model}, we introduce the model and recall the guarantee of greedy with estimated values. In Section~\ref{sec:rcm}, we prove Theorem~\ref{thm:rcm-main}. We show why expected sales are not submodular, introduce the fluid relaxation as a submodular surrogate, and bound its gap from expected sales through online matching. In Section~\ref{sec:mcc}, we prove Theorem~\ref{thm:greedy}. We show that expected sales are submodular under the Markov chain model, construct the evaluation oracle, and analyze the greedy algorithm. We discuss open problems in Section~\ref{sec:conclusion}.

\section{Model and Preliminaries}\label{sec:model}

\paragraph{Products, Inventory and Customers.}
A seller offers a set of substitutable products $\cN=\{1,\dots,n\}$, all sold at the same price, which we normalize to one. We denote the outside (no-purchase) option by $0$ and write $\cN^+:=\cN\cup\{0\}$. Before the selling season begins, the seller commits to an initial inventory vector $\vecc=(c_1,\dots,c_n)\in\mathbb{Z}_+^n$, where $c_i$ is the number of stocked units of product~$i$. There is no replenishment. The availability set $A(\vecc):=\{i\in\cN : c_i\ge 1\}$ consists of the products in stock. We write $\mathbf{e}_i$ for the $i$th unit vector and $\vecu\le \vecc$ for the componentwise order. A random number $T$ of customers arrive one at a time, indexed by $\tau\in[T]:=\{1,\dots,T\}$. The number $T$ is independent of the customers' choices, and its distribution is given by its finite support and the corresponding probabilities. Customers' preferences are independent and identically distributed draws from a discrete choice model $\phi$, where $\phi(i,S)$ is the probability that a customer who faces the availability set $S\subseteq\cN$ chooses option $i\in S\cup\{0\}$. Customer~$\tau$ faces $A(\vecc(\tau))$, where $\vecc(\tau)$ is the inventory vector upon her arrival (so $\vecc(1)=\vecc$). If she chooses a product $i$, she buys one unit of it and $\vecc(\tau+1)=\vecc(\tau)-\mathbf{e}_i$. If she chooses the outside option, $\vecc(\tau+1)=\vecc(\tau)$. We write
\[
1-\phi(0,S)=\sum_{i\in S}\phi(i,S)
\]
for the probability that a customer facing $S$ makes a purchase. In Remark~\ref{rem:nonstationary}, we extend our result for ranking-based models, with a smaller factor, to non-stationary arrival processes, where customers' preferences are independent but not identically distributed.
\paragraph{Ranking-Based Choice Model.}
In a ranking-based choice model (RCM), each customer draws a preference ranking $\sigma$ from a distribution $(p_\sigma)_{\sigma\in\boldsymbol{\Sigma}}$, where $\boldsymbol{\Sigma}$ is the set of all rankings of $\cN^+$. We write $\sigma(i)<\sigma(j)$ when option $i$ is preferred to option $j$. A customer facing $S$ chooses her most preferred option in $S\cup\{0\}$. The \emph{consideration set} of $\sigma$ is $A_\sigma:=\{i\in\cN:\sigma(i)<\sigma(0)\}$. A customer with ranking $\sigma$ buys if and only if $A_\sigma\cap S\ne\emptyset$, so
\[
1-\phi(0,S)=\sum_{\sigma\in\boldsymbol\Sigma}p_\sigma\mathbf 1\{A_\sigma\cap S\ne\emptyset\}.
\]
The RCM is equivalent to random utility models~\citep{block1959random,farias2013nonparametric}. The ranking distribution is given either as an explicit list or through an exact polynomial-time oracle for $\phi(0,S)$.

\paragraph{Markov Chain Choice Model.}
A Markov chain (MC) choice model is defined by an initial-state distribution $\bfmu=(\mu_i)_{i\in\cN^+}$ and a transition matrix $P=(P_{ik})_{i,k\in\cN^+}$~\citep{blanchet2016markov}. The outside option $0$ is absorbing, and we assume that from every state the chain reaches $0$ with probability one. Equivalently, $I-P_{\cN\cN}$ is invertible, where $P_{\cN\cN}$ is the restriction of $P$ to the products. A customer facing $S$ starts at state $i$ with probability $\mu_i$. If the current state is in $S\cup\{0\}$, she chooses it. Otherwise, she moves to state $k$ with probability $P_{ik}$. Choice probabilities can be computed by solving a linear system~\citep{blanchet2016markov}. The MC model is a special case of the RCM. Indeed, ordering the states by their first visit along a trajectory of the chain, and placing unvisited products after the outside option in a fixed order, yields a random ranking whose highest-ranked available option is the first available option visited by the chain.

\paragraph{Objective.}
For a deterministic horizon $t\in\mathbb{Z}_+$ and an inventory vector $\vecc$, let $R_t(\vecc)$ denote the expected number of units sold to $t$ customers, starting from $\vecc$. The expectation is over the customers' choices. Since all prices equal one, $R_t(\vecc)$ is also the expected revenue. For a capacity $K\ge1$, the seller's problem is
\begin{equation}\label{eq:problem}
\OPT:=\max_{\vecc\in\cC}\; F(\vecc),\qquad F(\vecc)\;:=\; \mathbb{E}_T\bigl[R_T(\vecc)\bigr],\qquad \cC:=\Bigl\{\vecc\in\mathbb{Z}_+^n : \sum_{i\in\cN} c_i\le K\Bigr\}.
\end{equation}
\paragraph{Submodularity.}
A set function $f:2^U\to\mathbb R$ is \emph{monotone} if $f(J)\le f(J')$ for all $J\subseteq J'\subseteq U$, and \emph{submodular} if $f(J\cup\{u\})-f(J)\ge f(J'\cup\{u\})-f(J')$ for all $J\subseteq J'\subseteq U$ and $u\in U\setminus J'$. We view an inventory vector as a set in an expanded ground set of units. Let the expanded ground set be $U:=\{(i,r):i\in\cN,\ 1\le r\le K\}$ and, for $J\subseteq U$, let $\vecc(J)$ be the inventory vector with $c_i(J):=|\{r:(i,r)\in J\}|$. Then, $\vecc(J)\in\cC$ if and only if $|J|\le K$. Under the capacity constraint, $K$ copies of each product suffice. For other downward-closed constraints, such as those in Remark~\ref{rem:constraints}, $t_{\max}$ copies of each product suffice, because no product sells more than $t_{\max}$ units. We say that a function $H:\mathbb Z_+^n\to\mathbb R$ is monotone and submodular if the set function $H(\vecc(\cdot)):2^U\to\mathbb R$ is, for every $K$.\footnote{This is equivalent to DR-submodularity, which is the property of diminishing returns for functions on multisets~\citep{kapralov2013online}.}
 Since units of a product are interchangeable, $H$ is submodular if and only if
\begin{equation}\label{eq:dr}
H(\vecc+\mathbf{e}_i+\mathbf{e}_j)-H(\vecc+\mathbf{e}_i)\;\le\; H(\vecc+\mathbf{e}_j)-H(\vecc)\qquad\text{for all }\vecc\in\mathbb{Z}_+^n\text{ and }i,j\in\cN,
\end{equation}
where $i$ and $j$ may be the same product. Adding one unit at a time, \eqref{eq:dr} shows that the marginal value of a unit of $i$ is nonincreasing in the inventory vector. We note that nonnegative mixtures preserve monotonicity and submodularity.

\paragraph{Greedy with Estimated Values.}
To maximize a monotone submodular set function $f$ over sets of at most $K$ elements, greedy starts from $J_0:=\emptyset$ and, for $\tau=0,\ldots,K-1$, adds to $J_\tau$ an element $u_\tau$ that maximizes an estimate of $f(J_\tau\cup\{u\})$. Since units of a product are interchangeable, on the set of units greedy simply adds a unit of a product $j$ that maximizes an estimate of $H(\vecc^{(\tau)}+\mathbf e_j)$, where $\vecc^{(\tau)}:=\vecc(J_\tau)$. If every estimate is within $\eta$ of the true value, then $f(J_\tau\cup\{u_\tau\})\ge\max_{u\in U}f(J_\tau\cup\{u\})-2\eta$. Therefore, the following guarantee applies with $\eta_\tau=2\eta$. With exact values, it is the classical guarantee of greedy~\citep{nemhauser1978}.
\begin{lemma}[\citealt{streeter2008online}, Theorem~3]\label{lem:approx-greedy}
Let $f:2^U\to\mathbb R_+$ be monotone and submodular, and let $f^\star:=\max_{J\subseteq U:\,|J|\le K}f(J)$. If, for every $\tau=0,\ldots,K-1$, the element $u_\tau$ added to $J_\tau$ satisfies $f(J_\tau\cup\{u_\tau\})\ge\max_{u\in U}f(J_\tau\cup\{u\})-\eta_\tau$ for some $\eta_\tau\ge0$, then $f(J_K)\ge(1-1/e)f^\star-\sum_{\tau=0}^{K-1}\eta_\tau$.
\end{lemma}
Theorem~3 of \citet{streeter2008online} is stated for schedules of actions with durations, and the lemma above is a special case where every action takes unit time.

\section{The Ranking-Based Choice Model}\label{sec:rcm}
In this section, we prove Theorem~\ref{thm:rcm-main}. Throughout, the choice model is ranking-based, and until the proof of the theorem (at the end of Section~\ref{sec:rcm-gap}) we fix a deterministic number $t$ of customers. We start by revisiting an example that demonstrates why the expected sales need not be submodular.

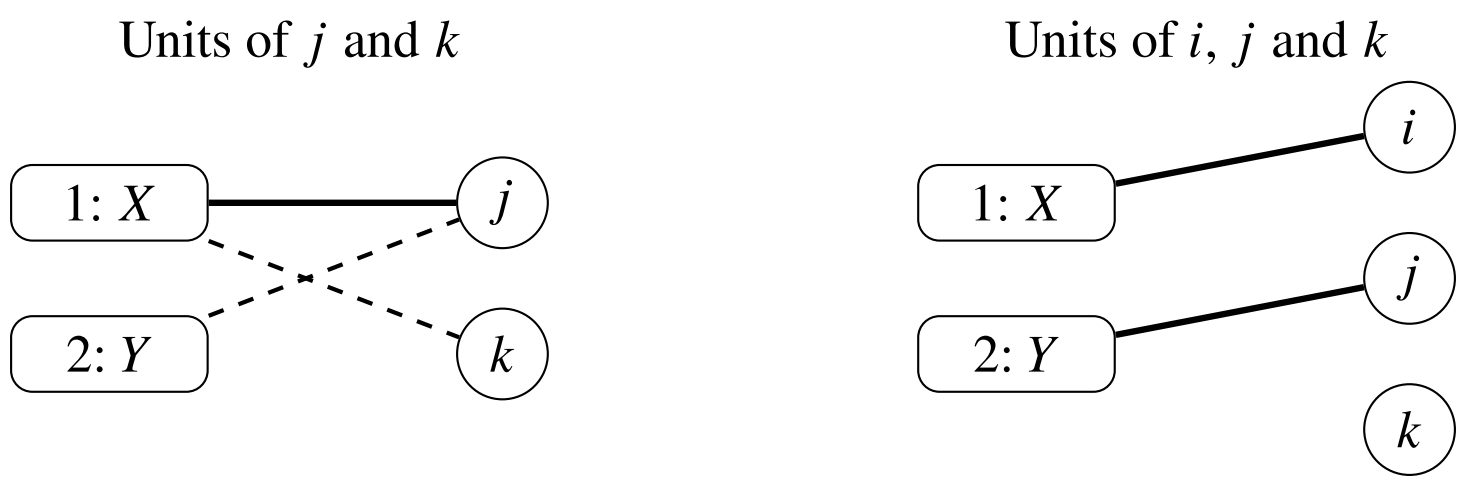
\begin{figure}[!ht]
\centering
\begin{tikzpicture}[font=\small,
  cust/.style={draw,rounded corners,minimum width=1.3cm,minimum height=0.5cm},
  unit/.style={draw,circle,minimum size=0.6cm,inner sep=0pt}]
\node at (1.2,1.55) {Units of $j$ and $k$};
\node[cust] (a1) at (0,0.5) {1: $X$};
\node[cust] (a2) at (0,-0.5) {2: $Y$};
\node[unit] (aj) at (2.6,0.5) {$j$};
\node[unit] (ak) at (2.6,-0.5) {$k$};
\draw[very thick] (a1) -- (aj);
\draw[dashed,thick] (a1) -- (ak);
\draw[dashed,thick] (a2) -- (aj);
\node at (7.2,1.55) {Units of $i$, $j$ and $k$};
\node[cust] (b1) at (6,0.5) {1: $X$};
\node[cust] (b2) at (6,-0.5) {2: $Y$};
\node[unit] (bi) at (8.6,1.0) {$i$};
\node[unit] (bj) at (8.6,0) {$j$};
\node[unit] (bk) at (8.6,-1.0) {$k$};
\draw[very thick] (b1) -- (bi);
\draw[very thick] (b2) -- (bj);
\end{tikzpicture}
\caption{The sample path of Example~\ref{ex:two-rankings} in which the first customer has ranking $X$ and the second has ranking $Y$, as a bipartite matching between customers (boxes) and units (circles). Thick lines show the greedy matching, that is, the purchases, and dashed lines show an optimal ex post matching. Left: with units of $j$ and $k$, the first customer buys $j$ and the second leaves, while the optimal matching serves both. Right: a unit of $i$ frees the unit of $j$ for the second customer, and the greedy matching is optimal.}
\label{fig:matching}
\end{figure}

\begin{example}[Two Rankings]\label{ex:two-rankings}
Consider a setting with three products $i$, $j$, $k$ and $t=2$ customers. Each customer independently draws one of two rankings: $X:=i\succ j\succ k\succ 0$ with probability $1/3$, and $Y:=j\succ 0$ with probability $2/3$.
\begin{itemize}
\item \emph{A unit of $i$ added to the set $\{k\}$} increases expected sales by $1/9$. It helps only when both customers have ranking $X$, in which case the first buys $i$ and the second buys $k$.
\item \emph{A unit of $i$ added to the set $\{j,k\}$ (one unit of each product)} increases expected sales by $2/9$ because it also helps when the first customer has ranking $X$ and the second has ranking $Y$ (Figure~\ref{fig:matching}). Without $i$, the first customer buys $j$ and the second leaves. With $i$, the first buys $i$ and the second buys $j$.
\end{itemize}
\end{example}

The same happens in the deterministic fluid version of stockout-based substitution, in which customers with rankings $X$ and $Y$ arrive as continuous flows at rates $1/3$ and $2/3$ over the time interval $[0,2]$ and buy their most preferred product in stock. There, the unit of $i$ increases sales by $0$ when added to $\{k\}$ and by $1/3$ when added to $\{j,k\}$. Using more sophisticated examples, one can show that the objective also does not satisfy weaker notions such as the submodular order property.

In Example~\ref{ex:two-rankings}, the unit of $i$ is worth more in the presence of $j$ because it frees the unit of $j$ for a later customer who wants only $j$. More generally, a unit can be consumed by customers with different rankings, and which of them consumes it is determined by the choice model rather than by the seller. Our main insight is that we can fix submodularity by \emph{controlling} the allocation of units to customers. In the example, if we could allocate units to customers after observing their rankings, then adding $i$ to $\{j,k\}$ would not increase sales at all, because we could give $k$ to the customer with ranking $X$ and $j$ to the customer with ranking $Y$ (the dashed lines in Figure~\ref{fig:matching}), whereas adding $i$ to $\{k\}$ would still increase sales by $1/9$. This simple observation suggests a natural surrogate function where we compute the optimal allocation of units to customers ex post. Suppose customer $\tau$ has ranking $\sigma$ and prefers the products in $A_\sigma$ to the outside option. Then the optimal ex post allocation is the solution to a maximum cardinality bipartite matching problem where on one side we have the units of the products and on the other we have the customers $\tau\in[t]$. Customer $\tau$ with ranking $\sigma$ has edges to the units of the products in $A_\sigma$. The optimal matching is thus a random variable, because the rankings (and, in general, the number of customers) are random. For simplicity of computing values, we define our surrogate as a well-known fluid LP upper bound on the expected value of the optimal matching.

The fluid LP is often used as a relaxation of the optimal offline value of the matching problem. As we will show, the actual expected sales are simply the expected size of the matching produced by a greedy algorithm, and existing results in the literature directly yield a tight approximation gap between the true objective and the surrogate that we define next.

\subsection{A Submodular Surrogate}\label{sec:rcm-fluid}
For any given $t$, the matching problem above is an instance of online bipartite matching with i.i.d.\ arrivals. On the offline side, product $i$ has $c_i$ units. On the online side, $t$ customers arrive one at a time, each with a ranking drawn independently from $(p_\sigma)_{\sigma\in\boldsymbol\Sigma}$, and a customer with ranking $\sigma$ can be matched to a unit of any product in $A_\sigma$. In expectation, $t\,p_\sigma$ of the $t$ customers have ranking $\sigma$. Let $E:=\{(i,\sigma)\in\cN\times\boldsymbol{\Sigma} : i\in A_\sigma\}$ be the set of edges between products and rankings. The \emph{fluid LP} is
\begin{equation}\label{eq:LP}
z_t(\vecc):=\max_{x\ge0}\biggl\{\sum_{(i,\sigma)\in E}t\,p_\sigma\,x_{i\sigma}:\ \sum_{i:(i,\sigma)\in E} x_{i\sigma}\le 1\ \ \forall\sigma\in\boldsymbol{\Sigma};\ \ \sum_{\sigma:(i,\sigma)\in E} t\,p_\sigma\,x_{i\sigma}\le c_i\ \ \forall i\in\cN\biggr\}.
\end{equation}
The variable $x_{i\sigma}$ is the probability that a customer with ranking $\sigma$ is matched to a unit of product $i$. These marginals sum to at most one for each ranking, and the expected number of customers matched to product $i$ is at most its inventory $c_i$. Every allocation policy induces marginals that satisfy these constraints, so the fluid LP upper bounds expected sales. This is a standard fact (see, for example, \citealt{golrezaei2014real}).

\begin{lemma}\label{lem:fluid-ub}
For every $\vecc\in\mathbb{Z}_+^n$ and every $t\in\mathbb{Z}_+$, $R_t(\vecc)\le z_t(\vecc)$.
\end{lemma}

Now, the fluid LP is, in fact, a transportation problem in which each unit of inventory is a source with one unit of supply and each ranking $\sigma$ is a sink with demand $t\,p_\sigma$. Its value is monotone in the set of units, because adding units relaxes the LP constraints. By a classical result of \citet{shapley1962complements}, it is submodular in the set of units.

\begin{lemma}[\citealt{shapley1962complements}; see \citealt{nemhauser1978}, Proposition~3.3]\label{lem:zlp-submod}
For every $t\in\mathbb{Z}_+$, $z_t$ is monotone and submodular.
\end{lemma}

Since nonnegative mixtures preserve monotonicity and submodularity, the surrogate
\(
g(\vecc):=\mathbb{E}_T[z_T(\vecc)]
\)
is monotone and submodular for any distribution of $T$, even though $F$ need not be submodular.

\subsection{The True Objective: Greedy Matching}\label{sec:rcm-gap}
For any instance (fixed graph), the actual objective is the expected value of an online greedy matching with deterministic tie-breaking. In particular, when a customer with ranking $\sigma$ arrives, consider the \emph{greedy} rule that matches her to an available unit of the product she ranks highest among the products in $A_\sigma$ that have units left, and leaves her unmatched if there is none. This is exactly stockout-based substitution: customers come one by one, and each chooses the first available product in her ranking.

The fluid LP is an upper bound, and how loose it can be is understood by comparing the worst-case performance of the online greedy policy with this offline benchmark. For any given $t$, \citet[Theorem~2.4]{devanur2019near} show, for a generalization of online bipartite matching with i.i.d.\ arrivals, that the greedy algorithm with arbitrary tie-breaking obtains in expectation at least a $1-1/e$ fraction of the optimal fractional solution of the expected instance, which here is the fluid LP. The factor $1-1/e$ cannot be improved (Remark~\ref{rem:tight}). This gives the following bound, which holds for every fixed horizon and is therefore independent of the distribution of $T$.

\begin{lemma}[\citealt{devanur2019near}, Theorem~2.4]\label{lem:stoch-gap}
For every $t\ge1$ and $\vecc\in\mathbb{Z}_+^n$, $R_t(\vecc)\ge(1-1/e)\,z_t(\vecc)$.
\end{lemma}

Taking expectations over $T$, which is independent of the customers' choices, gives
\begin{equation}\label{eq:gap-F}
F(\vecc)\ \ge\ \Bigl(1-\frac{1}{e}\Bigr)\, g(\vecc)\qquad\text{for every }\vecc\in\mathbb Z_+^n.
\end{equation}
Therefore, the surrogate approximates the true objective within a constant factor, and an inventory vector that maximizes the surrogate is a $(1-1/e)$-approximately optimal solution for the true objective.

\proof{Proof of Theorem~\ref{thm:rcm-main}.}
Suppose that we can compute $g(\vecc)$ exactly for any inventory vector $\vecc$. Then we can run greedy with exact values on the monotone submodular set function $g(\vecc(\cdot))$ on units (Lemma~\ref{lem:zlp-submod}). It returns $\hat\vecc\in\cC$, and Lemma~\ref{lem:approx-greedy} with $\eta_\tau=0$ gives
\begin{equation}\label{eq:up-sub}
g(\hat \vecc)
\ge(1-1/e)\max_{\vecc\in\cC}g(\vecc)
\ge(1-1/e)\OPT,
\end{equation}
where the second inequality uses $g\ge F$ (Lemma~\ref{lem:fluid-ub}). Then, by~\eqref{eq:gap-F}, $F(\hat \vecc)\ge(1-1/e)^2\OPT$.

It remains to show that we can evaluate $g$ efficiently. Greedy evaluates $g$ at $nK$ inventory vectors, and $g(\vecc)=\sum_{t\in\operatorname{supp}(T)}\mathbb P(T=t)\,z_t(\vecc)$, so greedy solves the fluid LP for $nK\,|\operatorname{supp}(T)|$ pairs of an inventory vector and a horizon. Hence, the running time is polynomial in $n$, $K$ and the size of the support of $T$, multiplied by the time needed to solve one fluid LP. Therefore, it suffices to show that we can compute $z_t(\vecc)$ (solve the fluid LP) efficiently. This is clear when the rankings with positive probability are listed explicitly, because the fluid LP then has polynomial size. When the choice model is specified through an oracle for $\phi(0,S)$, the probability of the outside option for any set $S$, and the rankings are not given explicitly (there may be $n!$ of them), we show that one can compute $z_t(\vecc)$ by submodular function minimization.

Specifically, we use the classical network-flow formulation of bipartite matching and the max-flow min-cut theorem~\citep[see, for example,][]{ahuja1993network}. Consider the network with a source, a node for each ranking $\sigma$, a node for each product $i$, and a sink. The network has an arc from the source to each ranking $\sigma$ with capacity $t\,p_\sigma$, an arc from $\sigma$ to each product in $A_\sigma$ with infinite capacity, and an arc from each product $i$ to the sink with capacity $c_i$. A flow that sends $t\,p_\sigma x_{i\sigma}$ units from $\sigma$ to $i$ corresponds to a feasible solution of the fluid LP, so the optimal value of the fluid LP is the maximum flow. If $S$ is the set of products on the sink side of a cut of finite capacity, every ranking $\sigma$ with $A_\sigma\cap S\neq\emptyset$ is also on the sink side. The total probability of these rankings is $\sum_{\sigma:\,A_\sigma\cap S\neq\emptyset}p_\sigma=1-\phi(0,S)$. Therefore, the max-flow min-cut theorem gives
\begin{equation}\label{eq:mincut}
z_t(\vecc)=\min_{S\subseteq\cN}\Bigl\{\sum_{i\notin S}c_i+t\,\bigl(1-\phi(0,S)\bigr)\Bigr\}.
\end{equation}
The function inside the minimum is submodular in $S$, because $1-\phi(0,S)=\sum_{\sigma}p_\sigma\mathbf 1\{A_\sigma\cap S\ne\emptyset\}$ is a weighted coverage function. Therefore,~\eqref{eq:mincut} is a submodular function minimization problem that needs only evaluations of $\phi(0,S)$. It can be solved with a number of oracle calls and other operations that is polynomial in $n$~\citep{iwata2001submodular}.\Halmos
\endproof

\begin{remark}[Other Constraints]\label{rem:rcm-ext}
The proof of Theorem~\ref{thm:rcm-main} uses only that $g$ is monotone and submodular and that $(1-1/e)\,g\le F\le g$. Therefore, any $\alpha$-approximation for maximizing $g$ under another constraint gives an $\alpha(1-1/e)$-approximation for the original problem. This gives the guarantees under budget and matroid constraints in Remark~\ref{rem:constraints}.
\end{remark}

\begin{remark}[Tightness]\label{rem:tight}
Each of the two factors $1-1/e$ is tight on its own. First, maximizing the surrogate is as hard as maximum coverage, so no polynomial-time algorithm can guarantee more than $1-1/e$ for it unless $\mathrm P=\mathrm{NP}$~\citep{feige1998threshold}. In particular, the fluid LP is equal to the expected sales when $T\equiv1$, i.e., $F(\vecc)= g(\vecc)=z_1(\vecc)$, and maximizing the expected sales when $T\equiv1$ is as hard as maximum coverage~\citep[Theorem~2]{goyal2016near}. Second, Lemma~\ref{lem:stoch-gap} is tight, for instance, when we have $t$ unit-capacity products and $t$ equally likely rankings that each place a different single product above the outside option. In this case, $z_t(\vecc)=t$, whereas $R_t(\vecc)$ is the expected number of distinct rankings among the $t$ customers, $t\bigl(1-(1-1/t)^t\bigr)$, which approaches $(1-1/e)\,t$. Whether both factors can be tight on the same instance is an intriguing question, which we leave open.
\end{remark}

\begin{remark}[Non-Stationary Arrival Process]\label{rem:nonstationary}
Suppose that customer $\tau$ draws her ranking from a distribution $(p^\tau_\sigma)_{\sigma\in\boldsymbol\Sigma}$ that may depend on $\tau$, independently of the other customers and of $T$. Then the guarantee of Theorem~\ref{thm:rcm-main} becomes $\frac12(1-1/e)\approx0.316$. With $t\,p_\sigma$ replaced by $\sum_{\tau\le t}p^\tau_\sigma$ in the fluid LP, Lemmas~\ref{lem:fluid-ub} and~\ref{lem:zlp-submod} still hold as is, and $g$ is still a monotone submodular upper bound on $F$. However, the true objective now corresponds to greedy for online bipartite matching with non-identical arrivals, and greedy is only $1/2$-competitive when arrivals are adversarial (arbitrary non-stationary)~\citep{karp1990, golrezaei2014real}.
\end{remark}

\begin{remark}[Independence from the Support of $T$]\label{rem:parametric}
The proof above solves the fluid LP once for every $t$ in the support of $T$. However, $g(\vecc)$ can be computed exactly with $O(n)$ submodular function minimizations, regardless of the support. By~\eqref{eq:mincut}, $z_t(\vecc)$ is a minimum of functions that are linear in $t$, so it is concave and piecewise linear in $t$. Since $1-\phi(0,S)$ is monotone in $S$, the largest minimizer in~\eqref{eq:mincut} can only shrink as $t$ increases~\citep{topkis1978minimizing}. In fact, the largest minimizers for the different values of $t$ form a nested chain $S_0\supsetneq S_1\supsetneq\cdots$ of subsets of $\cN$, and $z_t(\vecc)$ is linear in $t$ on each interval where the largest minimizer does not change. Therefore, there are at most $n+1$ pieces, and they can be found with $O(n)$ minimizations at rational values of $t$~\citep{eisner1976mathematical}. Since $z_t(\vecc)$ is linear in $t$ on each piece, computing $g(\vecc)$ requires only the probability and the partial mean of $T$ on at most $n+1$ intervals.
\end{remark}

\section{The Markov Chain Choice Model}\label{sec:mcc}

For the MC model, the probability $\phi(0,S)$ can be computed by solving a linear system, and thus Theorem~\ref{thm:rcm-main} directly applies and gives a $(1-1/e)^2$-approximation. However, this does not exploit the additional structure in the MC model where the customer rankings are not arbitrary but have a memoryless property. Recall that for MNL, a special case of the MC model, the true objective is known to be submodular for every distribution of $T$ using a coupling argument of \citet[Lemmas~3.6--3.8]{aouad2018greedy}, but the greedy algorithm could not be executed for MNL as is, because no efficient way to approximately evaluate the true objective was known until very recently~\citep{zhou2026optimal,zuo2026oneshot}.

We show more generally that the objective is monotone submodular under the MC model, generalizing the result for MNL (Section~\ref{sec:mc-submod}), and we develop an efficient way to approximately compute the expected sales of any given inventory vector (Section~\ref{sec:mult-oracle}). Together, these results give an improved $(1-1/e-\varepsilon)$-approximation for the MC model by implementing greedy directly on the true objective (end of Section~\ref{sec:mult-oracle}), which avoids the loss due to the surrogate objective and improves on $(1-1/e)^2$ for sufficiently small $\varepsilon$.

\subsection{Submodularity of Expected Sales}\label{sec:mc-submod}
Submodularity of expected sales is known in the static MC setting~\citep{desir2020constrained} and, under MNL, in the dynamic setting~\citep{aouad2018greedy}. The argument for MNL relies on the independence of irrelevant alternatives (IIA) property. Under MNL, a customer who finds her preferred product out of stock chooses like a new customer facing the remaining assortment. This does not extend directly to the MC model, because such a customer continues searching from that product in the Markov chain, and her subsequent choice probabilities need not coincide with those of a new arrival. Instead, our argument uses the memoryless property of the MC model, which implies that every customer continues from a given state in the same way. The main result of this subsection is the following.

\begin{theorem}\label{thm:mcc-submod}
For every integer $t\ge0$, the function $R_t$ is monotone and submodular:
\begin{equation}\label{eq:mc-dr-target}
R_t(\vecc+\mathbf{e}_i+\mathbf{e}_j)-R_t(\vecc+\mathbf{e}_i)\le R_t(\vecc+\mathbf{e}_j)-R_t(\vecc)
\qquad
(\vecc\in\mathbb{Z}_+^n,\ i,j\in\cN).
\end{equation}
\end{theorem}

Taking expectations over $T$, which is independent of the customers' choices, gives the following.

\begin{corollary}\label{cor:stoch-submod}
For every arrival distribution, $F(\vecc)=\mathbb{E}_T[R_T(\vecc)]$ is monotone and submodular.
\end{corollary}

We prove Theorem~\ref{thm:mcc-submod} by induction on the number of customers $t$. We couple the first customer's chain across the four inventory vectors $\vecc$, $\vecc+\mathbf e_i$, $\vecc+\mathbf e_j$ and $\vecc+\mathbf e_i+\mathbf e_j$ in~\eqref{eq:mc-dr-target}. Consider the first state that she reaches among the outside option and the products that are available in at least one of the four systems. If this state is the outside option, or a product that is available in all four systems, she makes the same choice in all four systems, and the induction hypothesis applies to the remaining $t-1$ customers. The difficult case is when she first reaches one of the added products, say $i$. Then, she buys it in the systems that have it and keeps searching in the others, where her later choices differ from those of a new customer. We handle this case using two properties of the MC model. First, one additional unit changes the inventory by at most one unit at all times (Lemma~\ref{lem:inventory-comparison}). Second, the order of the customers does not affect the distribution of the final inventory, so the customer who keeps searching can be moved to the end (Lemma~\ref{lem:reorder}). We prove these properties in Section~\ref{sec:mc-properties} and then the theorem in Section~\ref{sec:mc-proof}.
\subsubsection{Two Properties of the Markov Chain Model}\label{sec:mc-properties}

\paragraph{Notation.}
Throughout Sections~\ref{sec:mc-properties} and~\ref{sec:mc-proof}, the choice model is an MC model $(\bfmu,P)$, as defined in Section~\ref{sec:model}. A \emph{standard} customer starts her chain from $\bfmu$. Recall that $A(\vecu)$ is the set of products in stock under an inventory vector $\vecu$. For a state $x\in\cN^+$, let $\mathbb P(\vecu\prec0\mid x)$ denote the probability that a chain started at $x$ reaches a product in $A(\vecu)$ before $0$. In other words, $\mathbb P(\vecu\prec0\mid x)$ is the probability that a customer who starts at $x$ buys when the inventory is $\vecu$. It equals one if $x\in A(\vecu)$ and zero if $x=0$, and it is monotone in $\vecu$. For a horizon $t\ge1$, let $R_t(\vecu\mid x)$ denote the expected sales from $\vecu$ when the first of $t$ customers starts her chain at $x$ and the remaining $t-1$ customers are standard. Finally, let $\mathbf U_t(\vecu)$ denote the random inventory vector left after $t$ standard customers arrive and the starting inventory is $\vecu$.

The first property says that one additional unit changes the inventory by at most one unit at all times. It holds for any sequence of independent customers with transition matrix $P$, regardless of their initial-state distributions.

\begin{lemma}\label{lem:inventory-comparison}
For every $\vecu\in\mathbb Z_+^n$, $j\in\cN$, and $t\ge0$, the following hold.
\begin{enumerate}
\item $R_t(\vecu)\le R_t(\vecu+\mathbf e_j)\le R_t(\vecu)+1$.
\item For every state $x\in\cN^+$, $\E\bigl[\mathbb P\bigl(\mathbf U_t(\vecu+\mathbf e_j)\prec0\mid x\bigr)\bigr]\ge\E\bigl[\mathbb P\bigl(\mathbf U_t(\vecu)\prec0\mid x\bigr)\bigr]$.
\end{enumerate}
\end{lemma}

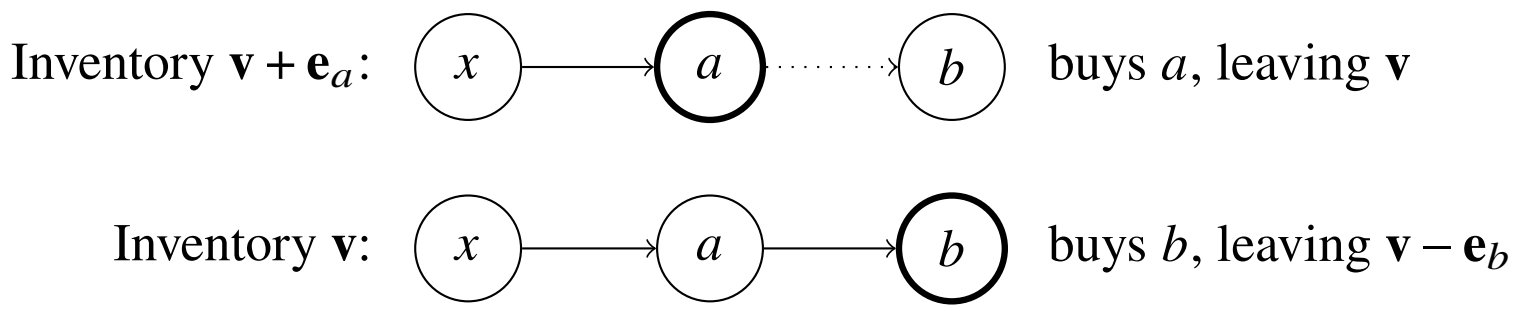
\begin{figure}[!ht]
\centering
\begin{tikzpicture}[font=\small,
  st/.style={draw,circle,minimum size=0.7cm,inner sep=0pt}]
\node[left] at (-0.5,0.6) {Inventory $\vecv+\mathbf e_a$:};
\node[st] (x1) at (0,0.6) {$x$};
\node[st,very thick] (a1) at (1.6,0.6) {$a$};
\node[st] (b1) at (3.2,0.6) {$b$};
\draw[->] (x1) -- (a1);
\draw[->,dotted] (a1) -- (b1);
\node[right] at (3.7,0.6) {buys $a$, leaving $\vecv$};
\node[left] at (-0.5,-0.6) {Inventory $\vecv$:};
\node[st] (x2) at (0,-0.6) {$x$};
\node[st] (a2) at (1.6,-0.6) {$a$};
\node[st,very thick] (b2) at (3.2,-0.6) {$b$};
\draw[->] (x2) -- (a2);
\draw[->] (a2) -- (b2);
\node[right] at (3.7,-0.6) {buys $b$, leaving $\vecv-\mathbf e_b$};
\end{tikzpicture}
\caption{Case~(b) of the proof of Lemma~\ref{lem:inventory-comparison}, when $v_a=0$ and the chain of the customer reaches $a$ before any other available product. Both systems follow the same chain path, and the thick circle marks the purchase. Afterwards, the two inventories again differ by one unit.}
\label{fig:coupling}
\end{figure}

\proof{Proof.}
Run the two systems starting from $\vecu+\mathbf e_j$ and $\vecu$ with the same customers, letting each customer follow the same path in the Markov chain in both systems. We show by induction over arrivals that the two inventories are either equal or of the form $\vecv+\mathbf{e}_a$ and $\vecv$, respectively, for some inventory $\vecv$ and product $a$. This holds initially. Identical inventories remain identical. Suppose the inventories before a customer arrives are $\vecv+\mathbf{e}_a$ and $\vecv$.
\emph{(a) $v_a>0$.} Both systems have the same available products, so the customer makes the same purchase in both and the difference remains $\mathbf{e}_a$.
\emph{(b) $v_a=0$.} The available sets differ only in $a$. If the chain reaches another available product or $0$ before reaching $a$, the customer makes the same choice in both systems. Otherwise, she buys $a$ in the first system, leaving $\vecv$, and continues her path after visiting $a$ in the second system (Figure~\ref{fig:coupling}). If she then buys a product $b$, the inventories become $\vecv$ and $\vecv-\mathbf e_b$. If she reaches $0$, both equal $\vecv$.

Therefore, on every sample path, the total sales of the two systems differ by $0$ or $1$. The total sales are monotone in the starting inventory, so this proves~(i). Moreover, the inventory of the first system is componentwise at least that of the second, so its available set is larger. Since $\mathbb P(\cdot\prec0\mid x)$ is monotone, this proves~(ii).\Halmos
\endproof

The second property says that the order of the customers does not matter for the distribution of the final inventory.

\begin{lemma}\label{lem:reorder}
Consider $m$ independent customers with transition matrix $P$ and arbitrary (possibly different) initial-state distributions. Given initial inventory $\vecu$, the distribution of the inventory left after all $m$ customers does not depend on the order of arrival. In particular, for every $x\in\cN^+$, $t\ge1$ and $\vecu\in\mathbb Z_+^n$,
\begin{equation}\label{eq:reorder}
R_t(\vecu\mid x)=R_{t-1}(\vecu)+\E\bigl[\mathbb P\bigl(\mathbf U_{t-1}(\vecu)\prec0\mid x\bigr)\bigr],
\end{equation}
where the left-hand side is the expected sales when the customer who starts at $x$ arrives first, and the right-hand side is the expected sales when she arrives last.
\end{lemma}

\proof{Proof.}
It suffices to show that exchanging two consecutive customers preserves the distribution of the inventory after the pair, because subsequent customers are independent of the pair, and their choices depend on earlier purchases only through the remaining inventory. Condition on the inventory $\vecv$ before the pair arrives and on the customers' initial states, and sample their paths until their first hits $x_1,x_2\in A(\vecv)\cup\{0\}$. If $x_1\ne x_2$, or one of them is $0$, the same products are bought in either order. If $x_1=x_2=a$ and $v_a\ge2$, both customers buy $a$ in either order. If $x_1=x_2=a$ and $v_a=1$, the customer served first buys the only unit of $a$, and the customer served second continues from $a$ facing $A(\vecv-\mathbf e_a)$. Because of the Markov property, this continuation has the same law regardless of which customer is served second. Therefore, under the natural coupling, both orders leave the same inventory after the pair.

For~\eqref{eq:reorder}, move the customer who starts at $x$ to the end. This does not change the distribution of the final inventory, and hence does not change expected sales. After this move, the other $t-1$ customers are standard and leave $\mathbf U_{t-1}(\vecu)$, and she then makes a purchase with probability $\mathbb P(\mathbf U_{t-1}(\vecu)\prec0\mid x)$.\Halmos
\endproof

\subsubsection{Proof of Theorem~\ref{thm:mcc-submod}}\label{sec:mc-proof}

Monotonicity follows from Lemma~\ref{lem:inventory-comparison}(i). We prove~\eqref{eq:mc-dr-target} by induction on $t$. It holds for $t=0$ because $R_0\equiv0$. Consider any $t\ge1$, $\vecc$ and $i,j\in\cN$, and suppose that $R_{t-1}$ satisfies~\eqref{eq:dr}.

We couple the first customer's choices under the four inventories $\vecc$, $\vecc+\mathbf{e}_i$, $\vecc+\mathbf{e}_j$ and $\vecc+\mathbf{e}_i+\mathbf{e}_j$ by sampling a single chain from $\bfmu$. Let $X$ be the first state the chain reaches in $A(\vecc+\mathbf{e}_i+\mathbf{e}_j)\cup\{0\}$. Before reaching $X$, the chain visits no product that is available in any of the four systems. Conditional on $X=x$, by the strong Markov property, the customer continues from $x$ as a fresh chain in every system, and the remaining $t-1$ customers are standard. Therefore, the conditional expected sales from each of the four inventories $\vecw$ equal $R_t(\vecw\mid x)$, and it suffices to show that, for every $x\in A(\vecc+\mathbf{e}_i+\mathbf{e}_j)\cup\{0\}$,
\begin{equation}\label{eq:mc-conditional-dr}
R_t(\vecc+\mathbf{e}_i+\mathbf{e}_j\mid x)-R_t(\vecc+\mathbf{e}_i\mid x) \le R_t(\vecc+\mathbf{e}_j\mid x)-R_t(\vecc\mid x).
\end{equation}
Since $A(\vecc+\mathbf{e}_i+\mathbf{e}_j)\cup\{0\}=A(\vecc)\cup\{0\}\cup\{i,j\}$, the following cases are exhaustive.

\emph{(a) $x\in A(\vecc)\cup\{0\}$.} If $x\in A(\vecc)$, the customer buys $x$ in all four systems. If $x=0$, she buys nothing in any of them. Let $\vecv:=\vecc-\mathbf e_x$ in the first case and $\vecv:=\vecc$ in the second. The first customer's sales cancel in~\eqref{eq:mc-conditional-dr}, which reduces to $R_{t-1}(\vecv+\mathbf{e}_i+\mathbf{e}_j)-R_{t-1}(\vecv+\mathbf{e}_i)\le R_{t-1}(\vecv+\mathbf{e}_j)-R_{t-1}(\vecv)$. This holds by the induction hypothesis.

\emph{(b) $x=i$ and $c_i=0$ (possibly $i=j$).} Under $\vecc+\mathbf{e}_i+\mathbf{e}_j$ and $\vecc+\mathbf{e}_i$, product $i$ is available, so the customer buys it, leaving $\vecc+\mathbf{e}_j$ and $\vecc$. Hence, the left-hand side of~\eqref{eq:mc-conditional-dr} equals $R_{t-1}(\vecc+\mathbf{e}_j)-R_{t-1}(\vecc)$. By Lemma~\ref{lem:reorder}, the right-hand side equals
\[
R_{t-1}(\vecc+\mathbf{e}_j)-R_{t-1}(\vecc)+\E\bigl[\mathbb P\bigl(\mathbf U_{t-1}(\vecc+\mathbf e_j)\prec0\mid i\bigr)\bigr]-\E\bigl[\mathbb P\bigl(\mathbf U_{t-1}(\vecc)\prec0\mid i\bigr)\bigr],
\]
and the difference of the last two terms is nonnegative by Lemma~\ref{lem:inventory-comparison}(ii).

\emph{(c) $x=j$, $c_j=0$ and $j\ne i$.} Inequality~\eqref{eq:mc-conditional-dr} is equivalent to $R_t(\vecc+\mathbf{e}_i+\mathbf{e}_j\mid x)-R_t(\vecc+\mathbf{e}_j\mid x) \le R_t(\vecc+\mathbf{e}_i\mid x)-R_t(\vecc\mid x)$, which is case~(b) with the roles of $i$ and $j$ exchanged.

Averaging~\eqref{eq:mc-conditional-dr} over $X$ completes the induction.\Halmos

\subsection{An Efficient Evaluation Oracle}\label{sec:mult-oracle}
As noted in prior work, finding an efficient way to evaluate the objective value or a close approximation appears to be challenging. Exact evaluation by dynamic programming may need to track up to $\prod_i(c_i+1)$ remaining-inventory states. The problem with direct Monte Carlo estimation is that the expected sales can be very small, which means that an accurate estimate requires a large number of samples. For example, if every customer buys with probability $p$, then a simulated sample path has no sale with probability at least $1-\E[T]\,p$, so a sample average needs about $1/(\E[T]\,p)$ samples merely to observe a sale, which is exponential in the input size when $p$ is exponentially small.

\begin{figure}[!ht]
\centering
\begin{tikzpicture}[font=\small,
  st/.style={draw,circle,minimum size=0.6cm,inner sep=0pt}]
\node at (1.6,2.0) {Chain of a customer facing $S$};
\draw[dashed,rounded corners] (2.0,-0.9) rectangle (3.6,1.3);
\node at (3.35,1.05) {$S$};
\node[st] (mu) at (0,0.2) {$\bfmu$};
\node[st] (s1) at (2.8,0.75) {$i$};
\node[st] (s2) at (2.8,-0.35) {$j$};
\node[st] (u1) at (1.1,1.2) {$k$};
\node[st] (z) at (1.1,-1.1) {$0$};
\draw[->] (mu) -- (u1);
\draw[->] (u1) -- (s1);
\draw[->] (mu) -- (s2);
\draw[->] (mu) -- (z);
\draw[->] (u1) to[bend left=10] (z);
\node at (7.4,2.0) {Conditioned on a sale};
\draw[dashed,rounded corners] (6.9,-0.9) rectangle (8.5,1.3);
\node at (8.25,1.05) {$S$};
\node[draw,rounded corners,inner sep=3pt] (mu2) at (4.6,0.2) {start};
\node[st] (t1) at (7.7,0.75) {$i$};
\node[st] (t2) at (7.7,-0.35) {$j$};
\node[st,draw=gray,text=gray] (z2) at (4.6,-1.1) {$0$};
\draw[->] (mu2) -- node[above,sloped,pos=0.4] {\scriptsize$\phi(i,S)/p_{\vecc}$} (t1);
\draw[->] (mu2) -- node[below,sloped,pos=0.4,yshift=-1pt] {\scriptsize$\phi(j,S)/p_{\vecc}$} (t2);
\end{tikzpicture}
\caption{Left: a customer facing the set $S=A(\vecc)$ of products in stock walks through the chain until she reaches $S$ or the outside option $0$. Coarsened to its two outcomes, the chain reaches $S$ with probability $p_{\vecc}=1-\phi(0,S)$. Right: conditional on reaching $S$, the customer starts directly at product $i\in S$ with probability $\phi(i,S)/p_{\vecc}$, and the customer never reaches $0$. Until the first sale, customers are independent trials with success probability $p_{\vecc}$. Therefore, we compute the probability of at least one sale exactly, and then simulate the conditioned process where the first buyer is guaranteed to visit $S$.}
\label{fig:coarse}
\vspace{-10pt}
\end{figure}
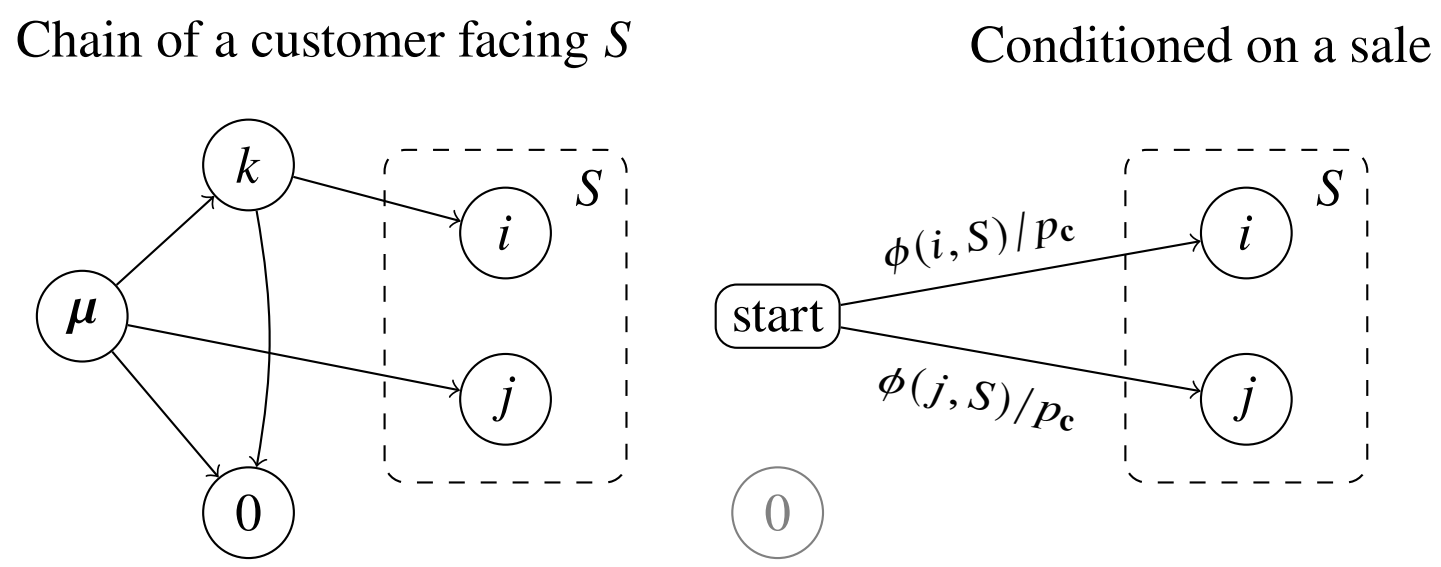

In other words, the trouble is that most simulated customers leave without buying and contribute nothing to the estimate. Our main idea therefore is quite simple: separate the outside option and simulate only the arrivals that lead to sales. For ease of exposition, we focus on the Markov chain choice model (we discuss generalizations in Remark~\ref{rem:oracle-general}). Here, we first consider a coarser chain where a customer facing $S=A(\vecc)$ has just two outcomes: she reaches $S$ and buys, or she reaches $0$ and leaves. Until the first sale occurs (after which something in $S$ might go out of stock), all customers reach $S$ independently with the same probability, so for any given horizon $t$, the probability that none of the $t$ customers reaches $S$ can be computed exactly. Now, we condition on the event that some customer reaches $S$ before $0$ and evaluate the number of sales on this event. The conditioned process is guaranteed to have at least one sale, which eliminates the main problem with direct sampling-based estimation (Figure~\ref{fig:coarse}). In other words, we can compute the probability of at least one sale exactly and simulate only conditional on it, and conditional on this event the number of sales lies between one and the total inventory.

In summary, the procedure has two steps. We compute the probability of at least one sale exactly (Section~\ref{sec:oracle-cond}) and sample the number of sales conditional on at least one sale (Section~\ref{sec:oracle-sim}). We bound the number of samples needed for an accurate estimate in Section~\ref{sec:oracle-acc}.

\vspace{-10pt}
\subsubsection{Conditioning on a Sale}\label{sec:oracle-cond}
Consider any inventory vector $\vecc$, and let $S:=A(\vecc)$, $C:=\|\vecc\|_1$, and $p_{\vecc}:=1-\phi(0,S)$. Let $Y(\vecc)$ denote the total number of units sold, so that $F(\vecc)=\mathbb E[Y(\vecc)]$. Until the first sale, the inventory does not change, so each customer independently makes a purchase with probability $p_{\vecc}$. Therefore, the probability of at least one sale is
\(
a_{\vecc}:=\mathbb P(Y(\vecc)\ge1)=\mathbb E_T\bigl[1-(1-p_{\vecc})^T\bigr],
\)
which can be computed exactly. Let $\zeta(\vecc):=\mathbb E[Y(\vecc)\mid Y(\vecc)\ge1]$ denote the expected number of sales conditional on at least one sale. Then, we have
\(
F(\vecc)=a_{\vecc}\,\zeta(\vecc)\) and \( 1\le\zeta(\vecc)\le C.
\)
If $a_{\vecc}=0$, then $F(\vecc)=0$ and the oracle returns $0$. Otherwise, it suffices to estimate $\zeta(\vecc)$.

\subsubsection{Conditional Simulation}\label{sec:oracle-sim}
To sample $Y(\vecc)$ conditional on $Y(\vecc)\ge1$, we simulate a complete sample path that begins with the first sale. The first step samples the first sale, and the second step simulates the remaining customers by repeating the same idea on the reduced inventory.
\begin{enumerate}
\item \emph{Sample the first sale.} Draw the horizon $t$ and the arrival $\kappa$ of the first sale with probability
\begin{equation}\label{eq:first-sale-law}
\mathbb P\bigl(T=t,\ \kappa=k\mid Y(\vecc)\ge1\bigr)=\frac{\mathbb P(T=t)\,(1-p_{\vecc})^{k-1}p_{\vecc}}{a_{\vecc}},\qquad t\in\operatorname{supp}(T),\ 1\le k\le t.
\end{equation}
Alternatively, we can first draw $t$ with probability $\mathbb P(T=t)[1-(1-p_{\vecc})^t]/a_{\vecc}$, and then draw $\kappa$ from the geometric distribution with parameter $p_{\vecc}$ truncated at $t$. Independently, draw the purchased product $i\in S$ with probability $\phi(i,S)/p_{\vecc}$.
\item \emph{Simulate the remaining customers.} Starting from inventory $\vecc-\mathbf e_i$, repeat this process for the remaining $t-\kappa$ customers. If $A$ is the current set of products in stock, the number of arrivals until the next sale is geometric with parameter $1-\phi(0,A)$, and the next customer who buys chooses product $i\in A$ with probability $\phi(i,A)/(1-\phi(0,A))$. The sample path ends when the horizon ends, the inventory runs out, or $\phi(0,A)=1$.
\end{enumerate}
The sample is the total number of sales on this path, which lies between one and $C$. Customers who make no purchase are never simulated individually, so each sample processes at most $C$ sales regardless of the length of the horizon. Moreover, the purchase probabilities need to be recomputed only when a product stocks out.

\subsubsection{Accuracy}\label{sec:oracle-acc}
The next lemma shows that an additive error of $\varepsilon\OPT$ can be achieved with a sample size independent of the magnitude of $\OPT$.

\begin{lemma}\label{lem:additive-f}
Consider any inventory vector $\vecc\in\mathbb Z_+^n$ with $a_{\vecc}>0$ and any $\varepsilon,\delta\in(0,1)$. Let $Y_1,\ldots,Y_N$ be independent samples generated by the conditional simulation, let $\hat\zeta(\vecc):=\frac1N\sum_{r=1}^N Y_r$, and let $\widehat F(\vecc):=a_{\vecc}\,\hat\zeta(\vecc)$. If $N\ge\bigl\lceil\frac{8C}{3\varepsilon^2}\log\frac2\delta\bigr\rceil$, then
\(
\mathbb P\bigl(|\widehat F(\vecc)-F(\vecc)|\le\varepsilon F(\vecc)\bigr)\ge1-\delta.
\)
In particular, if $\vecc\in\cC$, the estimate has additive error at most $\varepsilon\OPT$ with probability at least $1-\delta$.
\end{lemma}

\proof{Proof.}
We first show that each sample has the distribution of $Y(\vecc)$ conditional on $Y(\vecc)\ge1$. We start with the second step of the simulation. We claim that, for every inventory vector $\vecv$ and every number $s\ge0$ of customers, the second step started from $\vecv$ with $s$ customers generates the number of sales with the correct distribution. We prove the claim by induction on $\|\vecv\|_1$. If $\|\vecv\|_1=0$, there are no sales. Otherwise, let $A:=A(\vecv)$. The inventory remains $\vecv$ until the next sale, so the customers buy independently with probability $1-\phi(0,A)$. Therefore, the number of arrivals until the next sale is geometric with this parameter, and conditional on a purchase, the customer buys $i\in A$ with probability $\phi(i,A)/(1-\phi(0,A))$. After this sale, the remaining customers start from $\vecv-\mathbf e_i$, and the induction hypothesis applies.

For the first step, let $\kappa$ and $I$ denote the arrival and the product of the first sale in the original process. Using Bayes' rule, for every $t\in\operatorname{supp}(T)$, $1\le k\le t$ and $i\in S$, we have
\[
\mathbb P\bigl(T=t,\kappa=k,I=i\mid Y(\vecc)\ge1\bigr)=\frac{\mathbb P(T=t)(1-p_{\vecc})^{k-1}\phi(i,S)}{a_{\vecc}},
\]
which is the probability that the first step draws $t$, $k$ and $i$, by~\eqref{eq:first-sale-law}. Given the first sale, the remaining $t-k$ customers start from $\vecc-\mathbf e_i$, so the claim applies. Therefore, each $Y_r$ has the distribution of $Y(\vecc)$ conditional on $Y(\vecc)\ge1$. In particular, $\mathbb E[Y_r]=\zeta(\vecc)$ and $1\le Y_r\le C$.

Since $Y_r^2\le CY_r$, we have $\operatorname{Var}(Y_r)\le\mathbb E[Y_r^2]\le C\zeta(\vecc)$. Therefore, Bernstein's inequality gives
\begin{align}
\mathbb P\bigl(|\hat\zeta(\vecc)-\zeta(\vecc)|\ge\varepsilon\zeta(\vecc)\bigr)
&\le
2\exp\left(-\frac{N\varepsilon^2\zeta(\vecc)^2}{2C\zeta(\vecc)+(2/3)C\varepsilon\zeta(\vecc)}\right)
\le
2\exp\left(-\frac{3N\varepsilon^2}{8C}\right)
\le\delta,
\label{eq:conditional-estimator-concentration}
\end{align}
where the second inequality uses $\zeta(\vecc)\ge1$ and $\varepsilon<1$. Since $F(\vecc)=a_{\vecc}\zeta(\vecc)$ and $\widehat F(\vecc)=a_{\vecc}\hat\zeta(\vecc)$, with probability at least $1-\delta$ we have $|\widehat F(\vecc)-F(\vecc)|\le\varepsilon F(\vecc)$. For feasible $\vecc$, the inequality $F(\vecc)\le\OPT$ gives the additive guarantee.\Halmos
\endproof

We are now ready to prove Theorem~\ref{thm:greedy}. The proof applies a union bound to ensure that greedy has sufficiently accurate estimates throughout.

\proof{Proof of Theorem~\ref{thm:greedy}.}
Set $\varepsilon':=\varepsilon/(2K)$ and $\delta':=\delta/(nK)$. Greedy makes $nK$ calls to the procedure of Lemma~\ref{lem:additive-f}, with parameters $\varepsilon'$ and $\delta'$, each on an inventory vector with at most $K$ units, which is feasible. By Lemma~\ref{lem:additive-f}, each call returns an estimate $\widehat F(\vecu)$ of $F(\vecu)$ with $|\widehat F(\vecu)-F(\vecu)|\le\varepsilon'\OPT$ with probability at least $1-\delta'$. Since each evaluation uses independent simulations, the same bound holds conditional on previous evaluations and decisions. Therefore, a union bound shows that all estimates are accurate simultaneously with probability at least $1-nK\delta'=1-\delta$. On this event, Corollary~\ref{cor:stoch-submod} and Lemma~\ref{lem:approx-greedy}, applied to $F(\vecc(\cdot))$ with $\eta_\tau=2\varepsilon'\OPT=\varepsilon\OPT/K$, give
\[
F(\vecc^{(K)})\ge\Bigl(1-\frac1e\Bigr)\OPT-K\cdot\frac{\varepsilon\OPT}{K}=\Bigl(1-\frac1e-\varepsilon\Bigr)\OPT.
\]
For the running time, Lemma~\ref{lem:additive-f} requires $N=O\bigl(K(\varepsilon')^{-2}\log(1/\delta')\bigr)=O\bigl(K^3\varepsilon^{-2}\log(nK/\delta)\bigr)$ samples per call, so $O\bigl(nK^4\varepsilon^{-2}\log(nK/\delta)\bigr)$ samples in total. Each sample processes at most $K$ sales. The purchase probabilities are recomputed only when the available assortment changes, which occurs at most $n$ times. Finally, the conditional horizon and the geometric waiting times can be sampled without processing individual arrivals.\Halmos
\endproof
\begin{remark}[Other Choice Models]\label{rem:oracle-general}
	The procedure uses the choice model only through the choice probabilities $\phi(i,A)$ for the sets $A$ of products in stock. Moreover, the proof of Lemma~\ref{lem:additive-f} uses only that every sample sells between one and $C$ units. Therefore, the same procedure applies to any ranking-based model, or any other choice model, given an oracle for the choice probabilities $\phi(i,S)$. It estimates $F(\vecc)$ within a factor $1\pm\varepsilon$ with probability at least $1-\delta$ from the same number of samples. Each sample processes at most $C$ sales and queries the oracle for at most $n+1$ sets of products.
\end{remark}

\begin{remark}[Nonuniform Prices]\label{rem:prices} The idea also generalizes to nonuniform prices.
	Suppose that product $i$ has price $r_i\ge0$, and let $F_r$ and $\OPT_r$ denote the expected revenue and its maximum. Conditioning on a sale no longer controls the variance, because a high-priced product may be bought with small probability even conditional on some sale. Instead, we decompose revenue by product and condition, for each product, on an event whose probability does not depend on the inventory. Let $p_i:=\phi(i,\{i\})$. Since $\phi(i,A)\le p_i$ for every assortment $A$ that includes $i$, the customers' choices can be coupled with independent coins, one per customer, that come up heads with probability $p_i$, so that a customer buys $i$ only if her coin comes up heads. Let $E_i$ be the event that some coin comes up heads. Its probability $a_i=\E_T[1-(1-p_i)^T]$ is known. If $Y_i$ denotes the sales of $i$, then
	\(
	F_r(\vecc)=\sum_{i:c_i>0}v_i\,\E[Y_i/c_i\mid E_i]
	\)
	with known weights $v_i:=r_ic_ia_i$. Each conditional sell-through fraction lies in $[0,1]$ and can be sampled by conditional simulation. The first head follows a truncated geometric distribution. Before it, other products may sell with choice probabilities conditional on not buying $i$, and after it, customers follow the original choice process. Sampling a product with probability $v_i/V$, where $V:=\sum_{i:c_i>0}v_i$, and returning $V$ times its sell-through fraction gives an unbiased estimator of $F_r(\vecc)$ with values in $[0,V]$. Since $r_ia_i=F_r(\mathbf e_i)\le\OPT_r$, we have $V\le C\,\OPT_r$. Therefore, the variance is at most $V F_r(\vecc)\le C\,\OPT_r^2$, and Bernstein's inequality gives additive error $\varepsilon\OPT_r$ with probability at least $1-\delta$ from $O(C\varepsilon^{-2}\log(1/\delta))$ samples, as in Lemma~\ref{lem:additive-f}. Unlike Lemma~\ref{lem:additive-f}, the error is additive ($\varepsilon\OPT_r$) rather than multiplicative ($\varepsilon F_r(\vecc)$), because $p_i$ may overestimate the purchase probability of $i$, so $i$ may not sell even conditioning on $E_i$.
\end{remark}
\begin{remark}[Dependence on the Support of $T$]\label{rem:horizon}
The procedure of Section~\ref{sec:mult-oracle} sums over the support of $T$. By losing an $\varepsilon$ in the guarantee, we can replace $T$ by $\tilde T$, obtained by rounding $T$ down to a grid of $O(\varepsilon^{-1}\log t_{\max})$ points in which consecutive points differ by a factor of at most $1+\varepsilon$. The expected sales $R_t(\vecc)$ are concave in $t$, because the purchase probability of each customer can only decrease as products stock out. Since $R_0\equiv0$, the ratio $R_t(\vecc)/t$ is nonincreasing in $t$. Therefore, rounding $T$ down decreases the expected sales of every inventory vector by a factor of at most $1+\varepsilon$, and Theorem~\ref{thm:greedy} applied to $\tilde T$ gives a $(1-1/e-2\varepsilon)$-approximation. If $T$ has infinite support but belongs to a nice family such as geometric or Poisson, we can truncate the distribution at $t_{\max}$ such that $\E[T-t_{\max}\mid T\ge t_{\max}]\le\varepsilon\,t_{\max}$ and lose another $\varepsilon$. This decreases the expected sales of every inventory vector by a factor of at most $1+\varepsilon$.
\end{remark}

\section{Conclusion and Open Problems}\label{sec:conclusion}

We showed that dynamic assortment optimization under uniform prices admits constant-factor approximations for every distribution of the number of customers. Specifically, we gave a $(1-1/e)^2$-approximation under the ranking-based choice model and a $(1-1/e-\varepsilon)$-approximation under the Markov chain choice model. Several directions remain open.

\paragraph{The Approximability of the Ranking-Based Model.}
The best factor achievable in polynomial time lies between $(1-1/e)^2$ and $1-1/e$. Each of our two losses is tight on its own (Remark~\ref{rem:tight}). Whether both can be tight on the same instance seems (to us) non-trivial, and our guarantee may not be tight.

\paragraph{Nonuniform Prices.}
Under nonuniform prices, expected revenue need not be submodular even in the static Markov chain setting~\citep{desir2020constrained}, so our approach does not extend directly. Our evaluation oracle does extend (Remark~\ref{rem:prices}), but this does not by itself resolve the optimization problem.

\paragraph{Beyond the Markov Chain Model.}
Our proof of Theorem~\ref{thm:mcc-submod} uses only two properties. First, one additional unit changes the inventory by at most one unit at all times, which holds under every ranking-based model. Second, customers can be reordered without changing the distribution of the final inventory, which uses the memoryless property of the MC model. We note that the proof does not use that customers share the initial-state distribution. With the induction stated for every sequence of independent customers with transition matrix $P$, it shows that the theorem continues to hold when initial-state distributions differ across customers. Identifying other choice models with this reordering property is another interesting direction.

\section*{Acknowledgements and Disclosure of AI Use}
We thank Danny Segev for encouraging us to think about this problem and for thoughtful remarks on an early version of the results. In particular, it was thanks to Danny that we made the connection to submodular minimization in the proof of Theorem~\ref{thm:rcm-main}.

All mathematical ideas in this paper are the product of the authors. ChatGPT (OpenAI) and Claude (Anthropic) were used to check proofs, refine and simplify the arguments, and polish the manuscript. The authors take full responsibility for all contents of the manuscript.

\bibliographystyle{informs2014}
\bibliography{ref}
\end{document}